\documentclass{amsart}
\usepackage{amsthm}
\usepackage{amsfonts}
\usepackage{amssymb,latexsym}
\usepackage{amsbsy}
\usepackage{amsxtra}
\usepackage{hyperref}
\usepackage{amsmath}
\usepackage{mathabx}
\usepackage{mathrsfs}
\usepackage{graphicx}

\theoremstyle{plain}
\newtheorem{theorem}{Theorem}
\newtheorem{lemma}{Lemma}
\newtheorem{proposition}{Proposition}

\begin{document}
 \title{An Endpoint Version of Uniform Sobolev inequalities}
 \author{Tianyi Ren, Yakun Xi, Cheng Zhang}

 \begin{abstract}
We prove an endpoint version of the uniform Sobolev inequalities in Kenig-Ruiz-Sogge \cite{KRS}. It was known that strong type inequalities no longer hold at the endpoints; however, we show that restricted weak type inequalities hold there, which imply the earlier classical result by real interpolation. The key ingredient in our proof is a type of interpolation first introduced by Bourgain \cite{Bourgain}. We also prove restricted weak type Stein-Tomas restriction inequalities on some parts of the boundary of a pentagon, which completely characterizes the range of exponents for which the inequalities hold.

 \end{abstract}

\maketitle
\section{Introduction}
In this paper, we consider the second order constant coefficient differential operator
\[P(D)=Q(D)+\sum_{j=1}^na_j\frac{\partial}{\partial x_j}+b,\]
where $Q(\xi)$ denotes a nonsingular real quadratic form on $\mathbb{R}^n$, $n\ge3$ which, for some $2\le k\le n$, is given by \[Q(\xi)=-\xi_1^2-...-\xi_k^2+\xi_{k+1}^2+...+\xi_n^2,\]
$D=-i(\partial/\partial x_1,...,\partial/\partial x_n)$, and $a_1,...,a_n,b$ are complex numbers.
If $k=n$, the operator $P(D)$ has principal part $Q(D)=\Delta$ and is called elliptic. Otherwise, it is called non-elliptic.

Uniform Sobolev inequalities
\begin{equation}\label{sob}\|u\|_{L^{q}(\mathbb{R}^n)}\leqslant C\|P(D)u\|_{L^{p}(\mathbb{R}^n)},\quad u\in W^{2, p}(\mathbb{R}^n),\end{equation}
have been of interest to the study of unique continuation for partial differential equations. Here the constant $C$ should depend only on $n$ and $p$. If $P(D)=\Delta$, \eqref{sob} is just the classical Sobolev inequality. For more general elliptic operators, Kenig-Ruiz-Sogge \cite[Theorem 2.2]{KRS}  characterized the optimal range of exponents $(p,q)$ for which \eqref{sob} holds. Indeed, they showed that \eqref{sob} holds for elliptic $P(D)$ if and only if
$(p,q)$ satisfies the two conditions

(i)$\ \frac{1}{p}-\frac{1}{q}=\frac{2}{n},$

\begin{figure}
  \centering
    \includegraphics[height=7.5cm]{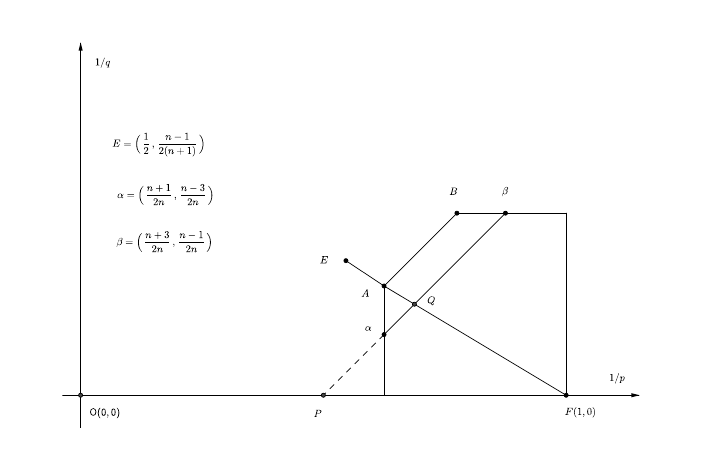}
\caption{The interpolation diagram for the resolvent estimates}
  \label{fig1}
\end{figure}

(ii)$\ \min\{|\frac{1}{p}-\frac{1}{2}|, |\frac{1}{q}-\frac{1}{2}|\}> \frac{1}{2n},$
\\
(i.e., $(1/p,1/q)$ lies on the open line segment joining $\alpha(\frac{n+1}{2n},\frac{n-3}{2n})$ and $\beta(\frac{n+1}{2n},\frac{n-1}{2n})$ in Figure \ref{fig1}). As pointed out in \cite[p.330]{KRS}, the chief technical difficulty in proving \eqref{sob} comes from the first order terms of $P(D)$. Indeed, \eqref{sob} follows from a localization argument and the uniform resolvent estimates:
\begin{equation}\label{model}
\|u\|_{L^{q}(\mathbb{R}^n)}\leqslant C\|(\Delta+z)u\|_{L^{p}(\mathbb{R}^n)}, \quad u\in W^{2, p}(\mathbb{R}^n),\ z\in \mathbb{C}.
\end{equation}

Since strong type inequality like \eqref{sob} and \eqref{model} no longer holds at the endpoints $\alpha$ and $\beta$, it is natural to ask whether restricted weak type inequality can be established at these endpoints. In this paper, we give a positive answer to this question when $P(D)$ is elliptic. 

\begin{theorem}
Let $n\geqslant3$. If $(1/p,1/q)=\alpha$ or $\beta$, then for any $z\in\mathbb{C}$, the inequality holds: \begin{equation} \label{15}
\|u\|_{L^{q,\infty}(\mathbb{R}^{n})}\leqslant C\|(\Delta+z)u\|_{L^{p, 1}(\mathbb{R}^n)},
\end{equation} where the constant $C$ depends only on $n$ and $p$.
\end{theorem}

\begin{theorem}
Let $n > 3$. If $(1/p,1/q)=\alpha$ or $\beta$, then there exists a constant $C$, depending only on $n$, such that whenever $P(D)$ is a second order constant coefficient differential operator with principal part $\Delta$, we have \begin{equation}\label{new}
\|u\|_{L^{q, \infty}(\mathbb{R}^{n})} \leqslant C \|P(D)u\|_{L^{p, 1}(\mathbb{R}^{n})}.
\end{equation}
\end{theorem}

A few remarks are in order. First, the above two theorems imply the corresponding classical results of Kenig-Ruiz-Sogge \cite{KRS} by real interpolation. Second, S. Guti\'errez \cite{Gutierrez} obtained restricted weak type resolvent estimates for the Laplacian at points $A(\frac{n+1}{2n}, \frac{(n-1)^2}{2n(n+1)})$ and $B(\frac{n^{2}+4n-1}{2n(n+1)},\frac{n-1}{2n})$ in Figure \ref{fig1}. The estimates cannot be uniform, but depend on $z$, because the exponent pairs are not on the line $\frac{1}{p}-\frac{1}{q}=\frac{2}{n}$. Finally, for non-elliptic $P(D)$, uniform restricted weak type estimates have been established in a recent work of Jeong-Kwon-Lee\cite{jkl}, which completely characterizes the range of $(p,q)$ for which \eqref{sob} holds in the non-elliptic case.

Theorem 2 follows from Theorem 1, a restricted weak type Stein-Tomas inequality (see Section 3) and the localization argument in \cite[p.335-337]{KRS}, after adapting several classical results to Lorentz spaces. To our knowledge, \eqref{new} is still open when $n=3$. The difficulty here is the failure of Littlewood-Paley inequality when an exponent becomes $\infty$ (see Proposition 1).

The rest of the paper is organized as follows. In Section 2, we prove Theorem 1 by an interpolation result first obtained by Bourgain \cite{Bourgain} and a variant of Stein's oscillatory integral theorem due to Sogge \cite{Osc}. The interpolation method of Bourgain was first brought to our attention by Bak-Seeger\cite{Seeger}. In Section 3, by a similar argument, we prove restricted weak type Stein-Tomas restriction inequalities on some parts of the boundary of a pentagon; this completely characterizes the range of exponents for which the inequalities hold. In Section 4, we prove Theorem 2 by establishing several classical results in the setting of Lorentz spaces and carrying out the localization argument in Kenig-Ruiz-Sogge \cite{KRS}.

\section{Proof of Theorem 1}

First of all, we need some reductions. It suffices to prove the theorem for one endpoint, say $\alpha(\frac{n+1}{2n}, \frac{n-3}{2n})$, because the other follows from duality. Furthermore, noting the gap condition $\frac{1}{p}-\frac{1}{q}=\frac{2}{n}$ on the exponents, we are able to reduce the theorem to the case where $z$ has unit length, $|z|=1$, after a simple rescaling argument. Finally, by continuity, we may assume that $\mathrm{Im}z\neq0$.

This last reduction enables us to study $(-|\xi|^{2}+z)^{-1}$, whose inverse Fourier transform is the fundamental solution of the operator ``$\Delta+z$'' in our theorem. Therefore, our theorem is a consequence of the following estimate for a multiplier operator: \begin{equation}
\Big\|\Big\{\frac{\hat{u}(\xi)}{-|\xi|^{2}+z}\Big\}^{\widecheck{\ }}\Big\|_{L^{\frac{2n}{n-3}, \infty}(\mathbb{R}^{n})}\leqslant C\|u\|_{L^{\frac{2n}{n+1}, 1}(\mathbb{R}^{n})}.
\end{equation}
This in turn, amounts to the inequality for a convolution operator \begin{equation}\label{16}
\Big\| u(x)\ast \Big\{\frac{1}{-|\xi|^{2}+z}\Big\}^{\widecheck{\ }}(x)\Big\|_{L^{\frac{2n}{n-3}, \infty}(\mathbb{R}^{n})}\leqslant C\|u\|_{L^{\frac{2n}{n+1}, 1}(\mathbb{R}^{n})}.\end{equation}
The proof of the theorem then transforms to the study of the kernel of this convolution operator: $K(x)=\Big\{\frac{1}{-|\xi|^{2}+z}\Big\}^{\widecheck{\ }}(x)$.
The expression for this kernel is actually already in literature, e.g. Gelfand-Shilov \cite{GF}:  \[
K(x)=\Big(\frac{z}{|x|^{2}}\Big)^\frac{n-2}{4}K_\frac{n-2}{2}(\sqrt{z|x|^{2}}),\] where \begin{equation} \label{10}
K_{\nu}(w)=\int_{0}^{\infty} e^{-w \mathrm{cosh}t} \mathrm{cosh}(\nu t) \mathrm{d}t
\end{equation} denotes the modified Bessel function. Along with the expression for $K(x)$, we will also need the following facts about the Bessel function, all of which are contained in
\cite[p.339]{KRS}.

First, a change of variable $u = e^{t}$ in the expression \eqref{10} for $K_{\nu}(w)$ immediately yields \begin{equation}\label{1}
|K_{\nu}(w)|\leqslant C|w|^{-|\mathrm{Re}(\nu)|},
\end{equation} for $|w|\leqslant1$ and $\mathrm{Re}(w)>0$, where the constant $C$ depends only on $\nu$.
Second, applying the formula (see \cite[p.19]{Erdelyi}) \begin{equation} \label{9}
\Gamma(\nu+\frac{1}{2})K_{\nu}(w) = \Big(\frac{\pi}{2w}\Big)^{\frac{1}{2}}e^{-w}\int_{0}^{\infty} e^{-t}t^{\nu-\frac{1}{2}}\Big(1+\frac{t}{2w}\Big)^{\nu-\frac{1}{2}} \mathrm{d}t,
\end{equation} which is valid when $\mathrm{Re}\nu \geqslant 0$, we obtain the behavior of the Bessel function for large $|w|$: \begin{equation}\label{2}
|K_{\nu}(w)|\leqslant Ce^{-\mathrm{Re}(w)}|w|^{-\frac{1}{2}}
\end{equation} whenever $|w|\geqslant1$ and $\mathrm{Re}(w)>0$.
Finally, formula \eqref{9} in fact tells us that \begin{equation}\label{3}
K_{\nu}(w)=a_{\nu}(w)w^{-\frac{1}{2}}e^{-w}
\end{equation} for $\mathrm{Re}(w)>0$, where the function $a_{\nu}(w)$ enjoys the decaying property \begin{equation} \label{4}
\Big| \Big( \frac{\partial}{\partial w}\Big)^{\alpha}a_{\nu}(w)\Big| \leqslant C_{\alpha}|w|^{-\alpha}.
\end{equation}

With these preparations, we embark on the task of proving the estimate \eqref{16} for the convolution operator $K(x)$. The idea is to treat the part of $K(x)$ inside the unit ball and the part outside separately, hence we break $K(x)=K_{1}(x)+K_{2}(x)$, where $K_{1}(x)$ is defined to equal $K(x)$ when $|x|\leqslant1$, and equal $0$ elsewhere.
By estimate \eqref{1}, considering the expression for $K(x)$, we easily obtain that \[
|K_{1}(x)|\leqslant C|x|^{-(n-2)}.\] Then the desired, indeed strong, inequality \begin{equation} \label{5}
\|u(x)\ast K_{1}(x)\|_{L^{\frac{2n}{n-3}}(\mathbb{R}^{n})}\leqslant C\|u\|_{L^{\frac{2n}{n+1}}(\mathbb{R}^{n})}
\end{equation} follows from Hardy-Littlewood-Sobolev inequality whenever the dimension $n>3$, noticing that the point $\alpha(\frac{n+1}{2n}, \frac{n-3}{2n})$ is on the line $\frac{1}{p}-\frac{1}{q}=\frac{2}{n}$. If $n=3$ however, we cannot apply Hardy-Littlewood Sobolev inequality, because one exponent $\frac{2n}{n-3}$ is $\infty$ then. Nevertheless, the restricted weak type estimate \begin{equation} \label{12}
\|u(x)\ast K_{1}(x)\|_{L^{\infty}(\mathbb{R}^{3})} \leqslant \|u\|_{L^{\frac{3}{2}, 1}(\mathbb{R}^{3})}\|K_{1}\|_{L^{3, \infty}(\mathbb{R}^{3})} \leqslant C\|u\|_{L^{\frac{3}{2}, 1}(\mathbb{R}^{3})}
\end{equation} still holds, by the Holder's inequality for Lorentz spaces. See for instance, \cite[Theorem 3.6]{O'Neil}.

After that, we turn to our analysis of $K_{2}(x)$, the part of $K(x)$ away from the origin. Applying \eqref{2} yields the estimate\[
|K_{2}(x)|\leqslant C|x|^{-\frac{n-1}{2}}e^{-|x|\mathrm{cos}(\frac{1}{2}\mathrm{arg}z)}.\]
Because of the exponential term, which may have the desired decaying property, we separate the case where $\mathrm{arg}z\in [-\frac{\pi}{2}, \frac{\pi}{2}]$, from the case where $\mathrm{arg}z\notin [-\frac{\pi}{2}, \frac{\pi}{2}]$. For the former situation, as just mentioned, the effect of the exponential decay yields the strong estimate\begin{equation} \label{6}
\|u(x)\ast K_{2}(x)\|_{L^{\frac{2n}{n-3}}(\mathbb{R}^{n})}\leqslant C\|u\|_{L^{\frac{2n}{n+1}}(\mathbb{R}^{n})},
\end{equation} which follows from Young's inequality.

The difficult situation is the latter one, and this is where Bourgain's interpolation comes into play. As pointed out before, this interpolation technique first appeared in \cite{Bourgain} when Bourgain was proving an endpoint bound for the spherical maximal function, and we first noticed it in \cite{Seeger}. There is also an elaboration on the abstract theory, developed for fairly general normed vector spaces, in Carbery-Seeger-Wainger-Wright \cite{Carbery}. We are not going into such abstractness and generality, but would rather state the result in our specific setting, that of $L^{p}$ spaces.

\begin{lemma}
Suppose that an operator $T$ between function spaces is the sum of the operators ${T_{j}}:$\[T=\displaystyle\sum_{j=1}^{\infty}T_{j}.\] If for $1 \leqslant p_{1}, p_{2}, q_{1}, q_{2} \leqslant \infty$, there exist $\beta_1,\beta_2 > 0$ and $M_{1}, M_{2} > 0$ such that each $T_{j}$ satisfies
\[\|T_j\|_{L^{p_1}\rightarrow L^{q_1}}\le M_12^{-j\beta_1},\]
and
\[\|T_j\|_{L^{p_2}\rightarrow L^{q_2}}\le M_22^{j\beta_2},\]
then we have restricted weak type estimate for the operator $T$ between two intermediate spaces:\[
\|Tf\|_{L^{q, \infty}}\leqslant C(\beta_{1}, \beta_{2}) M_{1}^{1-\theta}M_{2}^{\theta}\|f\|_{L^{p, 1}},\] where
\[\theta=\frac{\beta_{1}}{\beta_{1}+\beta_{2}},\]
\[\frac{1}{p}=\frac{1-\theta}{p_{1}}+\frac{\theta}{p_{2}}, \ \  \frac{1}{q}=\frac{1-\theta}{q_{1}}+\frac{\theta}{q_{2}}.\]
\end{lemma}

We return to the proof of the main theorem, dealing with the second situation where $\mathrm{arg}z\notin[-\frac{\pi}{2}, \frac{\pi}{2}]$. By the expression \eqref{3} for the Bessel function, \[
K_2(x)=|x|^{-\frac{n-1}{2}}e^{-i|x|\mathrm{sin}(\frac{1}{2}\mathrm{arg}z)}e^{i\frac{n-3}{4}\mathrm{arg}z}e^{-|x|\mathrm{cos}\frac{1}{2}\mathrm{arg}z}a_{\frac{n-2}{2}}(|x|e^{i\frac{1}{2}\mathrm{arg}z}),\]
where $a_{\frac{n-2}{2}}(w)$ satisfies the decaying property \eqref{4}.

We dyadically decompose the kernel $K_{2}(x)$. Fix a smooth function $\eta(x)$ that has support in $\{ x:|x|<1 \}$ and is equal to $1$ for $|x|<\frac{1}{2}$. Denote $\delta(x)=\eta(x)-\eta(2x)$. Then let $\beta_{0}(x)=\eta(x)$, and for each $j>1$, let $\beta_{j}(x)=\delta(2^{-j}x)$. It is easy to verify that $\sum_{j=0}^{\infty}\beta_{j}(x)=1$.

For each $j\geqslant0$, consider the operator $T_{j}$ given by the kernel \[ K_{2,j}(x)=\beta_{j}(x)K_{2}(x),\] i.e. $T_{j}u=u\ast K_{2, j}$. We need invoke the following variant of Stein's oscillatory integral theorem.

\begin{lemma}
Let $n\geqslant3$. Suppose that $1\leqslant p \leqslant 2$, $q=\frac{n+1}{n-1}p'$; in other words, the pair of exponents $(p, q)$ lies on the closed line segment joining $E(\frac{1}{2}, \frac{n-1}{2(n+1)})$ and $F(1, 0)$. Then, given a kernel of the form \[L(x)=\delta(x)b(x)e^{i\lambda|x|}|x|^{-\frac{n-1}{2}},\]where $\lambda \neq 0$, $\delta(x)$ is a smooth function supported in $\{x \in \mathbb{R}^{n}: \frac{1}{4} \leqslant |x| \leqslant 1\}$, $b(x)\in C^{\infty}(\mathbb{R}^{n})$, and $|((\frac{\partial}{\partial x})^{\alpha}b)(x)|\leqslant C_{\alpha}|x|^{-|\alpha|}$, we have the inequality \[
\|L\ast f\|_{L^{q}(\mathbb{R}^{n})}\leqslant C|\lambda|^{-\frac{n}{q}}\|f\|_{L^{p}(\mathbb{R}^{n})},\] where the constant $C$ depends only on the function $\delta(x)$ and finitely many of the $C_{\alpha}$ above.
\end{lemma}

This oscillatory integral theorem, proved by Sogge \cite{Osc}, follows from Stein's oscillatory integral theorem \cite{Stein}. See also \cite[p.341]{KRS}. It holds for pairs of exponents lying on the closed line segment $EF$ in Figure \ref{fig1}. Because of this, we are tempted to interpolate between the point $P(\frac{2}{n}, 0)$ on the $\frac{1}{p}$ axis and the point $Q(\frac{n^{2}+n+2}{2n^{2}}, \frac{n^{2}-3n+2}{2n^{2}})$, which is the intersection of the line $\alpha\beta$ and the line $EF$. See Figure \ref{fig1}. When $n=3$ however, $\alpha(\frac{n+1}{2n}, \frac{n-3}{2n})$ goes down to the $\frac{1}{r}$ axis and coincides with $P$, so interpolating between $P$ and $Q$ cannot produce a restricted weak type inequality at $\alpha$. Fortunately, we are able to remedy it by interpolating between two other points.

\textbf{Case $n>3$:} At $P$, since $K_{2,j}(x)\leqslant C|x|^{-\frac{n-1}{2}}$ by \eqref{2} and $K_{2,j}(x)$ is supported in $\{x:2^{j-2}\leqslant|x|\leqslant2^{j}\}$, Young's inequality yields \[\|K_{2,j}(x)\ast u(x)\|_{L^{\infty}(\mathbb{R}^{n})}
\leqslant C 2^{j\frac{n-3}{2}}\|u\|_{L^{\frac{n}{2}}(\mathbb{R}^{n})}.\]
At $Q$, we seek to prove \begin{equation} \label{11}
\|K_{2,j}(x)\ast u(x)\|_{L^{\frac{2n^{2}}{n^{2}-3n+2}}(\mathbb{R}^{n})} \leqslant C2^{-\frac{1}{n}j}\|u\|_{L^{\frac{2n^{2}}{n^{2}+n+2}}(\mathbb{R}^{n})}.
\end{equation} Changing the scale, replacing $x$ with $2^{j}x$, we would be done if we could show \[
\|\tilde{K}_{2,j}(x)\ast u(x)\|_{L^{\frac{2n^{2}}{n^{2}-3n+2}}(\mathbb{R}^{n})} \leqslant C2^{-\frac{n^{2}-3n+2}{2n}j}\|u\|_{L^{\frac{2n^{2}}{n^{2}+n+2}}(\mathbb{R}^{n})},\] where
\[\tilde{K}_{2,j}(x) = \delta(x)|x|^{-\frac{n-1}{2}}e^{-i2^{j}|x|\mathrm{sin}(\frac{1}{2}\mathrm{arg}z)}b(2^{j}x),\]
\[b(x) = e^{i\frac{n-3}{4}\mathrm{arg}z}e^{-|x|\mathrm{cos}\frac{1}{2}\mathrm{arg}z}a_{\frac{n-2}{2}}(|x|e^{i\frac{1}{2}\mathrm{arg}z}),\]
and $\delta(x)$ is as mentioned in the dyadic decomposition. The kernel $\tilde{K}_{2,j}$ is easily seen to fall within the hypotheses of the above oscillatory integral theorem, remembering that $a_{\frac{n-2}{2}}(w)$ satisfies the decaying property \eqref{4}. Hence we obtain the inequality \eqref{11} by applying Lemma 2, with the constant $C$ independent of $j$.

The estimates at $P$ and $Q$ for the operator $T_{j}$ enables us to utilize Bourgain's interpolation, resulting in the desired restricted weak type estimate for the operator $T$, which is the sum of the $T_{j}$. Specifically,
\[\theta=\frac{\frac{1}{n}}{\frac{1}{n}+\frac{n-3}{2}}=\frac{2}{n^{2}-3n+2};\]
\[\Big(\frac{n^{2}+n+2}{2n^{2}}, \frac{n^{2}-3n+2}{2n^{2}}\Big)\cdot(1-\theta)+\Big(\frac{2}{n}, 0\Big)\cdot\theta=\Big(\frac{n+1}{2n}, \frac{n-3}{2n}\Big).\]
Noticing that the last pair of exponents is precisely the endpoint $\alpha$, we have by Lemma 1
\begin{equation} \label{7}
\|K_{2}(x)\ast u(x)\|_{L^{\frac{2n}{n-3}, \infty}(\mathbb{R}^{n})}\leqslant C\|u\|_{L^{\frac{2n}{n+1}, 1}(\mathbb{R}^{n})},
\end{equation} where, independent of $z\in\mathbb{C}$, the constant $C$ depends exclusively on the dimension $n$. \eqref{7} together with \eqref{5} and \eqref{6} gives the conclusion of Theorem 1 whenever $n>3$.

\textbf{Case $n=3$:} We interpolate instead between the two points $O(0, 0)$ and $F(1, 0)$, bearing in mind that the oscillatory integral theorem Lemma 2 holds for the latter pair too. At $O$, Young's inequality gives \[
\|K_{2,j}(x)\ast u(x)\|_{L^{\infty}(\mathbb{R}^{n})}\le C2^{2j}\|u\|_{L^{\infty}(\mathbb{R}^{n})},\]
while a change of scale argument as in the above case along with the oscillatory integral theorem Lemma 2 shows immediately that at $F$,\[
\|K_{2,j}(x)\ast u(x)\|_{L^{\infty}(\mathbb{R}^{n})}
\leqslant C2^{-j}\|u\|_{L^{1}(\mathbb{R}^{n})}.\] Then we compute
\[\theta=\frac{1}{1+2}=\frac{1}{3},\]
\[(0, 0)\cdot\theta+(1, 0)\cdot(1-\theta)=(\tfrac{2}{3}, 0).\]
Again, the last pair of exponents is our target pair. That concludes our proof of the Theorem 1. \qed

\section{Endpoint Version of Stein-Tomas Fourier Restriction Theorem}

\begin{figure}
  \centering
    \includegraphics[height=7.5cm]{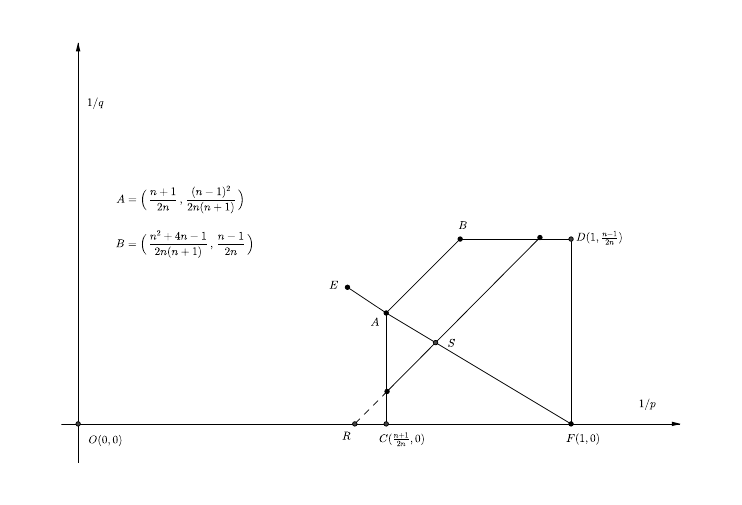}
\caption{The interpolation diagram for the restriction estimates}
  \label{fig2}
\end{figure}

The original Stein-Tomas restriction theorem \cite{Stein}, \cite{Tomas} states that if the dimension $n\geqslant3$, one has the inequality \begin{equation}\label{strong}
\Big\|\int_{S^{n-1}}\hat{f}(\xi)\mathrm{e}^{2\pi i\langle x, \xi\rangle}\,\mathrm{d}\sigma(\xi)\Big\|_{L^{q}(\mathbb{R}^{n})}\leqslant C\|f\|_{L^{p}(\mathbb{R}^{n})},
\end{equation}at the point $(1/p,1/q)=(\frac{n+3}{2n+2},\frac{n-1}{2n+2})$, which is the midpoint of $AB$ in Figure \ref{fig2}. Sogge \cite{Osc} extended this result in showing that the same inequality holds for pairs of exponents $(p, q)$ off the line of duality satisfying $1 \leqslant p<\frac{2n}{n+1}, q=\frac{n+1}{n-1}p'$. They constitute the half open line segment connecting $F(1, 0) $ and $A(\frac{n+1}{2n}, \frac{(n-1)^{2}}{2n(n+1)})$ (exclusive). Therefore by duality and interpolation, the Stein-Tomas restriction inequality is true for pairs of exponents in the interior of the pentagon in Figure \ref{fig2}. Furthermore, Bak-Seeger \cite[Proposition 2.1]{Seeger}) established restricted weak type inequalities
\begin{equation}\label{weak}
\Big\|\int_{S^{n-1}} \hat{f}(\xi)\mathrm{e}^{2\pi i \langle x, \xi \rangle}\mathrm{d}\sigma(\xi)\Big\|_{L^{q,\infty}(\mathbb{R}^{n})}\leqslant C\|f\|_{L^{p,1}(\mathbb{R}^{n})}
\end{equation} at the vertices $A$ and $B$. Then by real interpolation, strong type inequalities as \eqref{strong} hold on the open segment $AB$. In addition, strong type inequalities trivially hold on the half open segments $CF$ and $DF$ (excluding $C$ and $D$) by Young's inequality.

However, no results seem to have been established on $AC$ and $BD$ before. It is clear that strong Stein-Tomas can not hold on these two segments. Indeed, in \cite{Vilenkin}, radial functions belonging to $L^{\frac{2n}{n+1}}{(\mathbb{R}^{n})}$ are constructed that have infinite Fourier transforms on $S^{n-1}$. Moreover, neither strong type nor restricted weak type inequality holds outside of the pentagon in Figure \ref{fig2}. In fact, if there were a restricted weak type inequality somewhere outside this pentagon, then by real interpolation, we would either get a strong inequality on the line of duality $q=p'$ with $p > \frac{2n+2}{n+3}$, or get a strong inequality somewhere on $AC$ or $BD$. This is a contradiction, remembering that the range $1 \leqslant p \leqslant \frac{2n+2}{n+3}$ is sharp for a strong restriction estimate on the line of duality (see \cite[p.387, 2.1.1]{HA}).

In this section, we show that restricted weak type inequality as \eqref{weak} holds on the closed segments $AC$ and $BD$, by an argument similar to the proof of Theorem 1. With this result and the discussion above, we completely characterize the range of $(p,q)$ for which either strong Stein-Tomas or restricted weak type Stein-Tomas holds.

\begin{theorem}
Let $n \geqslant 3$. If $p=\frac{2n}{n+1}$ and $\frac{2n(n+1)}{(n-1)^{2}} \leqslant q \leqslant \infty$, or $1\leqslant p\leqslant \frac{2n(n+1)}{n^{2}+4n-1}$ and $q=\frac{2n}{n-1}$, then
\[\Big\|\int_{S^{n-1}}\hat{f}(\xi)\mathrm{e}^{2\pi i\langle x, \xi\rangle}\,\mathrm{d}\sigma(\xi)\Big\|_{L^{q, \infty}(\mathbb{R}^{n})}\leqslant C\|f\|_{L^{p, 1}(\mathbb{R}^{n})}\]
\end{theorem}

\begin{proof}
Our result follows from an analysis of the convolution operator whose kernel is the Fourier transform of the Lebesgue measure on the unit sphere, like in the classical case. This kernel is well-known to have the expression
\[
K(x)=2\pi |x|^{-\frac{n-2}{2}}J_{\frac{n-2}{2}}(2\pi|x|),\] where $J_{\nu}(w)$ is the Bessel function, see for instance, \cite[p.347-348]{HA}. Later, we will need the following fact about $J_{\nu}(w)$ for $\nu = m$ positive, integral or half integral, and $w = r$ real, positive and greater than $1$, which is also well-known: it takes the form \begin{equation} \label{8}
J_{m}(r) = \sum_{\pm} r^{-\frac{1}{2}}e^{\pm ir}a_{\pm}(r),
\end{equation} where the functions $a_{\pm}(r)$, $r > 1$ are smooth and satisfy the decay property
\[\Big|\frac{\mathrm{d}^{k}}{\mathrm{d}r^{k}}a_{\pm}(r)\Big| \leqslant C_{k}r^{-k}.\]
This expression can be found in \cite[p.338]{HA}; see also \cite[Theorem 1.2.1]{fio}. Again, dyadically decompose $K(x)$, letting $\beta_{j}(x)$, $j \geqslant 0$ be as in the decomposition before and $K_{j}(x)=\beta_{j}(x) K(x)$.

We first treat the case of the vertex at $C$, because it is exceptional. For this, we wish to apply Bourgain's interpolation to the point $O(0, 0)$ and the point $F(1, 0)$. There is no need to worry about the part $K_{0}(x)$ of $K(x)$ near the origin, since $K(x)$ is the Fourier transform of a compactly supported distribution and is thus smooth. Away from the origin, i.e., for $j > 0$, at $O$, Young's inequality gives \[
\|K_{j}(x)\ast f(x)\|_{L^{\infty}(\mathbb{R}^{n})}\leqslant C2^{\frac{n+1}{2}j}\|f\|_{L^{\infty}(\mathbb{R}^{n})},\]
while at $F$, still applying Young's inequality \[
\|K_{j}(x)\ast f(x)\|_{L^{\infty}(\mathbb{R}^{n})} \leqslant C2^{-\frac{n-1}{2}j}\|f\|_{L^{1}(\mathbb{R}^{n})}.\]
With the following interpolation computation,
\[\theta=\frac{\frac{n-1}{2}}{\frac{n-1}{2}+\frac{n+1}{2}}=\frac{n-1}{2n},\]
\[(1, 0)\cdot (1-\theta)+(0, 0)\cdot \theta=(\tfrac{n+1}{2n}, 0),\]
we obtain the restricted weak type Stein-Tomas inequality at $C (\frac{n+1}{2n}, 0)$: \[
\Big\|\int_{S^{n-1}} \hat{f}(\xi) \mathrm{e}^{2\pi i\langle x, \xi\rangle} \mathrm{d}\sigma(\xi)\Big\|_{L^{\infty}(\mathbb{R}^{n})}\leqslant C\|f\|_{L^{\frac{2n}{n+1}, 1}(\mathbb{R}^{n})}.\] Duality then produces the same restricted weak type inequality for the pair of exponents $D (1, \frac{n-1}{2n})$. However here, we cannot apply real interpolation to the points $A$ and $C$, nor can we apply it to $B$ and $D$, since they are along a vertical or horizontal line, which violates a hypothesis of the real interpolation theorem.

Nevertheless, we can proceed as in the proof of Theorem 1 to obtain a restricted weak type Stein-Tomas inequality at every point on the line segment joining $A$ and $C$ and its dual line segment joining $B$ and $D$. Indeed, for each $\frac{2}{n+1} < k <\frac{n+1}{2n}$, we interpolate between the point $R(k, 0)$ and the point $S=(\frac{n-1}{2n}+\frac{n+1}{2n}k, \frac{n-1}{2n}(1-k))$, which is the intersection of the line $\frac{1}{p}-\frac{1}{q}=k$ and the line $EF$. At $R$ for each $j > 0$, \[
\|K_{j}(x)\ast f(x)\|_{L^{\infty}(\mathbb{R}^{n})} \leqslant C2^{j(n(1-k)-\frac{n-1}{2})}\|f\|_{L^{\frac{1}{k}}(\mathbb{R}^{n})},
\] while at $S$ for $j > 0$, the familiar change of scale argument in which we replace $x$ with $2^{j}x$, together with Lemma 2 produces \[
\|K_{j}(x)\ast f(x)\|_{L^{s}(\mathbb{R}^{n})} \leqslant C2^{j(-\frac{n+1}{2}k+1)}\|f\|_{L^{r}(\mathbb{R}^{n})},\] where $r$ and $s$ denote the exponents corresponding to $S$. Finally we verify
\[\theta = \frac{\frac{n+1}{2}k-1}{\frac{n-1}{2}(1-k)},\]
\[\Big(\frac{n-1}{2n}+\frac{n+1}{2n}k, \frac{n-1}{2n}(1-k)\Big)(1-\theta)+(k, 0)\theta = \Big(\frac{n+1}{2n}, 1-k-\frac{n-1}{2n}\Big),\]
which gives us a restricted weak type inequality \eqref{weak} at every point on the line segment $AC$, as we hoped. Duality then produces the same results for the dual line segment $BD$.
\end{proof}

\section{Proof of Theorem 2}

As noted in the introduction, the procedure is essentially the same as in \cite{KRS}, with a restricted weak type Stein-Tomas inequality and the adaptations of several classical results to Lorentz spaces.

\paragraph{\textbf{1. Littlewood-Paley inequality}}
For $t \in \mathbb{R}$, let $\chi(t)$ be the characteristic function of the set $\{t: |t| \in [1, 2]\}$, and let $\chi_{k}(\xi_{n}) = \chi(2^{k}\xi_{n})$. Suppose $g$ is any function in $\mathscr{S}(\mathbb{R}^{n})$ for simplicity.
\begin{proposition}
For any $1 < p < \infty$ and $1 \leqslant q \leqslant \infty$, there exist constants $C_1, C_2$, depending only on $p$, $q$ and $n$, such that the inequalities below hold
\[C_{2}\|g\|_{L^{p, q}(\mathbb{R}^{n})} \leqslant \Big\|\Big( \sum_{k=-\infty}^{\infty} |\{\chi_{k}(\xi_{n})\hat{g}(\xi)\}^{\widecheck{\ }}|^{2}\Big)^{\frac{1}{2}} \Big\|_{L^{p, q}(\mathbb{R}^{n})} \leqslant C_{1}\|g\|_{L^{p, q}(\mathbb{R}^{n})}.\]
\end{proposition}

\begin{proof} The upper bound follows directly from the usual Littlewood-Paley inequality and real interpolation. The lower bound can be obtained by imitating the duality argument in \cite[p.105]{Stein2}.
\end{proof}

\paragraph{\textbf{2. Minkowski's Inequality}}

\begin{proposition}
\[
\Big\| \Big(\sum_{k=-\infty}^{\infty} |F_{k}(x)|^{2} \Big)^{\frac{1}{2}}\Big\|_{L^{s, \infty}(\mathbb{R}^{n})} \leqslant C\Big(\sum_{k=-\infty}^{\infty} \|F_{k}(x)\|_{L^{s, \infty}(\mathbb{R}^{n})}^{2} \Big)^{\frac{1}{2}},\] for any $s >2$, where the $C$ depends only on $s$ and $n$;
\[\Big(\sum_{k=-\infty}^{\infty} \|F_{k}(x)\|_{L^{r, 1}(\mathbb{R}^{n})}^{2} \Big)^{\frac{1}{2}} \leqslant C\Big\| \Big(\sum_{k=-\infty}^{\infty} |F_{k}(x)|^{2} \Big)^{\frac{1}{2}}\Big\|_{L^{r, 1}(\mathbb{R}^{n})},\] for any $1 < r < 2$, where the $C$ depends only on $r$ and $n$.
\end{proposition}

\begin{proof} The first inequality follows easily from Minkowski's inequality, recalling that $L^{p, q}$ is a Banach space when $1<p<\infty$, $1\leqslant q\leqslant\infty$. The second inequality results from a standard duality argument.
\end{proof}

\paragraph{\textbf{3. Stein-Tomas Inequality}}

\begin{proposition}
If $n \geqslant 3$, then
\[\Big\|\int_{S^{n-1}}\hat{f}(\xi)\mathrm{e}^{2\pi i\langle x, \xi\rangle}\,\mathrm{d}\sigma(\xi)\Big\|_{L^{\frac{2n}{n-3}, \infty}(\mathbb{R}^{n})}\leqslant C\|f\|_{L^{\frac{2n}{n+1}, 1}(\mathbb{R}^{n})}\]
\end{proposition}

\begin{proof}This is a special case of Theorem 3.
\end{proof}

\paragraph{\textbf{4. H\"ormander's Multipliers Theorem}}

\begin{proposition}
Suppose that $m \in L^{\infty}(\mathbb{R}^{n})$ satisfies, for some integer $s > \frac{n}{2}$,
\[\sum_{0 \leqslant |\alpha| \leqslant s} \sup_{\lambda > 0} \lambda^{-n}\|\lambda^{|\alpha|}D^{\alpha}\beta(\cdot/\lambda)m(\cdot)\|_{L^{2}(\mathbb{R}^{n})}^{2} < \infty,\]
whenever $\beta \in C_{0}^{\infty}(\mathbb{R}^{n} \setminus \{0\})$. Then for $1 < p < \infty$ and $1 \leqslant q \leqslant \infty$, the inequality holds
\[\|T_{m}f\|_{L^{p, q}(\mathbb{R}^{n})} \leqslant C_{p, q}\|f\|_{L^{p, q}(\mathbb{R}^{n})},\]
where the $T_{m}$ is the multiplier operator with multiplier $m(x)$:
\[T_{m}f = \{m(\xi)\hat{f}(\xi)\}^{\widecheck{\ }}.\]
\end{proposition}

\begin{proof}This is a consequence of H\"ormander's multiplier theorem and real interpolation.
\end{proof}

Having the above four propositions at hand, we are able to prove Theorem 2, following the localization argument in \cite[p.335-337]{KRS}. Below is a sketch of this argument.

By a reduction process similar to that at the beginning of the proof of Theorem 1 (see \cite[p.335]{KRS}), it suffices to prove for Theorem 2 the following special case \begin{equation}
\|u\|_{L^{\frac{2n}{n=3}, \infty}(\mathbb{R}^{n})} \leqslant C \Big\| \Big\{\Delta+1+\epsilon \Big(\frac{\partial}{\partial x_{n}}+i\beta \Big) \Big\}u \Big\|_{L^{\frac{2n}{n+1}, 1}(\mathbb{R}^{n})},
\end{equation} where $\epsilon \neq 0$, $\beta \neq 0$. Theorem 2 is then a consequence of the estimate for a multiplier operator below: \begin{equation} \label{13}
\Big\|\Big\{\frac{\hat{f}(\xi)}{-|\xi|^{2}+1+i\epsilon(\xi_{n}+\beta)}\Big\}^{\widecheck{\ }}\Big\|_{L^{\frac{2n}{n-3}, \infty}(\mathbb{R}^{n})} \leqslant C\|f\|_{L^{\frac{2n}{n+1}, 1}(\mathbb{R}^{n})}.
\end{equation}
Denote the multiplier in \eqref{13} by $m(\xi)$. Also, for $t \in \mathbb{R}$, let $\chi(t)$ be the characteristic function of the set $\{t: |t| \in [1, 2]\}$, and set $\chi_{k}(\xi_{n}) = \chi(2^{k}(\xi_{n}+\beta)$. For convenience, denote $\chi_{k}(\xi_{n})m(\xi)$ as $m_{k}(\xi)$. It then suffices to prove a similar estimate for the multiplier $m_{k}(\xi)$: \begin{equation} \label{14}
\|\{m_{k}(\xi)\hat{f}(\xi)\}^{\widecheck{\ }}\|_{L^{\frac{2n}{n-3}, \infty}(\mathbb{R}^{n})} \leqslant C\|f\|_{L^{\frac{2n}{n+1}, 1}(\mathbb{R}^{n})}.
\end{equation} Indeed, noting that $\frac{2n}{n+1} < 2 < \frac{2n}{n-3}$, if we had estimate \eqref{14}, we may apply the second part of Proposition 1, the first part of Proposition 2, the second part of Proposition 2, and the first part of Proposition 1, in that order, to obtain
\[\begin{aligned}
\|\{m(\xi)\hat{f}(\xi)\}^{\widecheck{\ }}\|_{L^{\frac{2n}{n-3}, \infty}(\mathbb{R}^{n})} & \leqslant C\Big\|\Big(\sum_{k=-\infty}^{\infty
} |\{m_{k}(\xi)\hat{f}(\xi)\}^{\widecheck{\ }}|^{2}\Big)^{\frac{1}{2}}\Big\|_{L^{\frac{2n}{n-3}, \infty}(\mathbb{R}^{n})}\\
& \leqslant C\Big(\sum_{k=-\infty}^{\infty} \|\{m_{k}(\xi)\hat{f}(\xi)\}^{\widecheck{\ }}\|_{L^{\frac{2n}{n-3}, \infty}(\mathbb{R}^{n})}^{2}\Big)^{\frac{1}{2}}\\
& \leqslant C\Big(\sum_{k=-\infty}^{\infty} \|\{\chi_{k}(\xi_{n})\hat{f}(\xi)\}^{\widecheck{\ }}\|_{L^{\frac{2n}{n+1}, 1}(\mathbb{R}^{n})}^{2}\Big)^{\frac{1}{2}}\\
& \leqslant C\Big\|\Big(\sum_{k=-\infty}^{\infty
} |\{\chi_{k}(\xi_{n})\hat{f}(\xi)\}^{\widecheck{\ }}|^{2}\Big)^{\frac{1}{2}}\Big\|_{L^{\frac{2n}{n+1}, 1}(\mathbb{R}^{n})}\\
& \leqslant C\|f\|_{L^{\frac{2n}{n+1}, 1}(\mathbb{R}^{n})},\\ \end{aligned} \] which is the result we are seeking.

To prove inequality \eqref{14}, we first apply the special case we just proved, \eqref{15} in Theorem 1 to $z = 1+i\epsilon2^{-k}$ and obtain \begin{equation}
\Big\|\Big\{\frac{\chi_{k}(\xi_{n})\hat{f}(\xi)}{-|\xi|^{2}+1+i\epsilon2^{-k}}\Big\}^{\widecheck{\ }}\Big\|_{L^{\frac{2n}{n-3}, \infty}(\mathbb{R}^{n})} \leqslant C\|f\|_{L^{\frac{2n}{n+1}, 1}(\mathbb{R}^{n})}.
\end{equation} By taking difference, it then remains only to demonstrate the inequality \begin{equation}
\begin{split}
&\Big\|\Big\{\frac{\chi_{k}(\xi_{n})[i\epsilon(\xi_{n}+\beta-2^{-k})]\hat{f}(\xi)}{(-|\xi|^{2}+1+i\epsilon(\xi_{n}+\beta))(-|\xi|^{2}+1+i\epsilon2^{-k})}\Big\}^{\widecheck{\ }}\Big\|_{L^{\frac{2n}{n-3}, \infty}(\mathbb{R}^{n})} \\
& \leqslant C\|f\|_{L^{\frac{2n}{n+1}, 1}(\mathbb{R}^{n})}.\\
\end{split}
\end{equation}
Now if we use polar coordinates $\xi=\rho\omega$, we will get, after applying Minkowski's inequality for the Lorentz space $L^{\frac{2n}{n-3}}(\mathbb{R}^{n})$ and Proposition 3, the following string of inequalities \begin{align}
{}& \Big\|\Big\{\frac{\chi_{k}(\xi_{n})[i\epsilon(\xi_{n}+\beta-2^{-k})]\hat{f}(\xi)}{(-|\xi|^{2}+1+i\epsilon(\xi_{n}+\beta))(-|\xi|^{2}+1+i\epsilon2^{-k})}\Big\}^{\widecheck{\ }}\Big\|_{L^{\frac{2n}{n-3}, \infty}(\mathbb{R}^{n})} \notag \\
\begin{split}
{}& \leqslant \int_{0}^{\infty} \Big\|\int_{S^{n-1}} \frac{\epsilon\hat{f}(\rho\omega)\chi_{k}(\xi_{n})(\xi_{n}+\beta-2^{-k})e^{i\rho \langle \omega, x \rangle}}{(-\rho^{2}+1+i\epsilon(\xi_{n}+\beta))(-\rho^{2}+1+i\epsilon2^{-k})}\mathrm{d}\omega\Big\|_{L^{\frac{2n}{n-3}, \infty}(\mathbb{R}^{n})}\rho^{n-1}\mathrm{d}\rho \notag
\end{split} \\
\begin{split}
{}& \leqslant C\int_{0}^{\infty} \rho\Big\|\Big\{\frac{\epsilon\hat{f}(\xi)\chi_{k}(\xi_{n})(\xi_{n}+\beta-2^{-k})}{(-\rho^{2}+1+i\epsilon(\xi_{n}+\beta))(-\rho^{2}+1+i\epsilon2^{-k})}\Big\}^{\widecheck{\ }}\Big\|_{L^{\frac{2n}{n+1}, 1}(\mathbb{R}^{n})}\mathrm{d}\rho. \notag
\end{split}
\end{align}
Finally, by Proposition 4, this last expression is majorized by
\[ C\|f\|_{L^{\frac{2n}{n+1}, 1}(\mathbb{R}^{n})} \int_{0}^{\infty} \frac{|\epsilon2^{-k}\rho|}{(\rho^{2}-1)^{2}+(\epsilon2^{-k})^{2}}\mathrm{d}\rho,\]
which is dominated by $C\|f\|_{L^{\frac{2n}{n+1}, 1}(\mathbb{R}^{n})}$. This concludes the proof of Theorem 2. \qed

One thing worthy of mentioning is that in Theorem 2, we exclude the case $n=3$, in contrast with the corresponding result in Kenig-Ruiz-Sogge \cite{KRS}. This is due to the fact that the exponent $\frac{2n}{n-3}$ is $\infty$ when $n=3$ and thus Littlewood-Paley inequality fails.

\vspace{0.5cm}

\textbf{Acknowledgement}. The authors would like to express their gratitude to their advisor, Professor Christopher D. Sogge, for bringing this research topic and Bourgain's interpolation to their attention, and also for the invaluable guidance and suggestions he provided.

\bibliography{Endpoint}

\bibliographystyle{plain}

\vspace{0.5cm}

\setlength{\parindent}{0cm}
\textit{Email address}: tyren@math.jhu.edu

\setlength{\parindent}{0cm}
Department of Mathematics, Johns Hopkins University, Baltimore, MD 21218, USA

\vspace{0.5cm}

\setlength{\parindent}{0cm}
\textit{Email address}: yxi4@math.rochester.edu

\setlength{\parindent}{0cm}
Department of Mathematics, University of Rochester, Rochester, NY 14627, USA

\vspace{0.5cm}

\setlength{\parindent}{0cm}
\textit{Email address}: czhang67@math.jhu.edu

\setlength{\parindent}{0cm}
Department of Mathematics, Johns Hopkins University, Baltimore, MD 21218, USA

\end{document}